\documentclass[10pt]{article}
\usepackage{amssymb, amsmath, url, graphicx, setspace, geometry}
\usepackage{calrsfs}
\usepackage{wasysym}

\def\3{\subset }
\def\4{\subseteq }
\def\<{\left<}
\def\>{\right>}

\def\bit{\begin{itemize}}
\def\eit{\end{itemize}}
\def\3{\subset }
\def\4{\subseteq }

\def\0{\leqno}

\def\barr{\begin{array}}
\def\earr{\end{array}}

\def\Z{{\rlap{$\kern2pt{\rm Z}$}{\rm Z}\,}}

\title{\bf Probabilistic aspects of ZM-groups}
\author{Mihai-Silviu Lazorec}
\date{December 17, 2017}

\begin{document}

\maketitle

\begin{abstract}
In this paper we study probabilistic aspects such as (cyclic) subgroup commutativity degree and (cyclic) factorization number of ZM-groups. We show that these quantities can be computed using the sizes of the conjugacy classes of these groups. In the end, we point out a class of groups whose (cyclic) subgroup commutativity degree vanishes asymptotically.
\end{abstract}



\section{Introduction}

The probabilistic aspects of finite groups theory constitute a topic which is regularly studied in the last years. The relevant results obtained in \cite{3} and \cite{5}-\cite{10} using the \textit{commutativity degree} of a finite group $G$, which measures the probability that two elements of the group commute, had a major impact on the research related to the probabilistic aspects associated to groups. Therefore, inspired by the previously mentioned concept, the \textit{subgroup commutativity degree} and the \textit{cyclic subgroup commutativity degree}  of $G$ were introduced in \cite{16} and \cite{19}, respectively. These quantities measure the probability that $G$ is an Iwasawa group, i.e. a nilpotent modular group, and they are defined by 
$$sd(G)=\frac{1}{|L(G)|^2}|\lbrace (H, K) \in L(G)^2 \ | \ HK=KH\rbrace|$$ and $$csd(G)=\frac{1}{|L_1(G)|^2}|\lbrace (H, K) \in L_1(G)^2 \ | \ HK=KH\rbrace|,$$ where $L(G)$ and $L_1(G)$ denote the subgroup lattice and the poset of cyclic subgroups of $G$, respectively. In other words, the (cyclic) subgroup commutativity degree represents the probability that two (cyclic) subgroups of $G$ commute. There is a major interest in finding explicit formulas for these two concepts under the assumption that $G$ belongs to a well known class of finite groups. As a consequence, one can study the asymptotic behavior of $sd(G)$ and $csd(G)$. For such kind of results, we refer the reader to \cite{16, 19, 21}. 

Another remarkable quantity associated to a finite group $G$ is its \textit{factorization number} which is denoted by $F_2(G)$. This number is obtained by counting the pairs $(H,K)\in L(G)^2$ satisfying $G=HK$. These pairs of subgroups are called factorizations of $G$. The factorization number of $G$ was studied in \cite{2,11,17} and it is strongly connected with $sd(G)$. More exactly, $$sd(G)=\frac{1}{|L(G)|^2}\sum\limits_{H\in L(G)}F_2(H).$$ 
In \cite{20}, the \textit{cyclic factorization number} of a finite group $G$, denoted by $CF_2(G)$, was introduced and it was obtained by counting the pairs $(H,K)\in L_1(G)^2$ such that $HK=G$. Obviously, the connection between $CF_2(G)$ and $csd(G)$ is similar to the one between $F_2(G)$ and $sd(G)$, and it is given by $$csd(G)=\frac{1}{|L_1(G)|^2}\sum\limits_{H\in L(G)}CF_2(H).$$
Our aim is to study these four probabilistic aspects for ZM-groups.

The paper is organized as follows. Section 2 recalls the structure of ZM-groups as well as some of the properties of this remarkable class of groups. Explicit formulas for the (cyclic) factorization number of ZM-groups are obtained in Section 3. Section 4 covers the main results of the paper. More exactly, we find an explicit formula that allows us to compute the (cyclic) subgroup commutativity degree of ZM-groups, we show that our results generalize the ones obtained for dihedral groups in \cite{16,19}, and we indicate a class of groups whose (cyclic) subgroup commutativity degree vanishes asymptotically. We end the paper by suggesting some open problems in the last section.

Most of our notation is standard and will usually not be repeated here. Elementary notions and results on groups can be found in \cite{4,13}. For subgroup lattice concepts we refer the reader to \cite{12,14,15}.  

\section{Properties of ZM-groups}

The aim of this Section is to recall some of the properties of ZM-groups that will be further used in our study. We will start by describing the structure of these groups, we continue by indicating the lattice of subgroups and the poset of cyclic subgroups for this type of groups and we discuss about its conjugacy classes. 

According to \cite{4}, a ZM-group is a finite group with all Sylow subgroups cyclic and it has the following structure $$ZM(m,n,r)=\langle a,b \ | \ a^m=b^n=1, \ b^{-1}ab=a^r\rangle,$$ where the parameters $(m,n,r)\in \mathbb{N}^{3}$ are satisfying the conditions $(m,n)=(m,r-1)=1$ and $r^n \equiv 1 \ (mod \ m).$
The order of $ZM(m,n,r)$ is $mn$ and the above properties of the triple $(m,n,r)$ force $m$ to be an odd positive integer. In \cite{1}, the subgroups of a ZM-group were completely described by showing that there is a bijection between the set 
$$L=\bigg\lbrace (m_1,n_1,s) \in \mathbb{N}^3 \ | \ m_1|m, \ n_1|n, \ 0\leq s<m_1, \ m_1|s\frac{r^n-1}{r^{n_1}-1} \bigg\rbrace$$ and $L(ZM(m,n,r))$. Moreover, for each triple $(m_1,n_1,s)\in L$, the corresponding subgroup of $ZM(m,n,r)$ is $$H(m_1,n_1,s)=\bigcup\limits_{i=1}^{\frac{n}{n_1}}\alpha (n_1,s)^{i}\langle a^{m_1}\rangle=\langle a^{m_1},\alpha(n_1,s)\rangle,$$ where $\alpha(x,y)=b^{x}a^{y}$, for $x\in\lbrace 0,1,2,\ldots, n-1\rbrace, y\in \lbrace 0,1,2,\ldots, m-1\rbrace$. Since $(\frac{m}{m_1}, \frac{n}{n_1})=1$, it follows that the order of a subgroup $H(m_1,n_1,s)\in L(ZM(m,n,r))$ is $\frac{mn}{m_1n_1}$. Also, according to \cite{18}, there is a bijection between the set $$L_1=\bigg\lbrace (m_1,n_1,s) \in L \ | \ \frac{m}{m_1}|r^{n_1}-1\bigg\rbrace$$ and the poset of cyclic subgroups $L_1(ZM(m,n,r))$. Hence, the (cyclic) subgroups of a ZM-group are fully described.

The conjugacy classes of a ZM-group will play an important role in our study. Therefore, we continue our short review on this class of groups, by recalling some properties of the set of conjugacy classes of an arbitrary group. We mention that the reader may find the proofs of the following facts in \textbf{Section III.3} of \cite{14}. We will denote the set of conjugacy classes of a group $G$ and the conjugacy class of a subgroup $H\in L(G)$ by $\mathcal{C}(G)$ and $[H]$, respectively. If $G$ is finite, then $(\mathcal{C}(G), \leq )$ is a poset, the partial ordering relation being defined as follows: for $[H_1], [H_2]\in \mathcal{C}(G)$, we have $[H_1]\leq [H_2]$ iff there exists $g \in G$ such that $H_1\subseteq H_2^{g}$. Moreover, a finite group $G$ has all Sylow subgroups cyclic if and only if $\mathcal{C}(G)$ is isomorphic to the lattice of all divisors of $|G|$. Also, two subgroups of a finite group $G$ whose Sylow subgroups are cyclic, are conjugate iff they have the same order. Applying these results for a ZM-group, we infer that $$\mathcal{C}(ZM(m,n,r))=\lbrace [H(m_i,n_i,0)] \ | \ m_i|m, n_i|n\rbrace.$$
Also, since $(m,n)=1$ and there is an isomorphism between $\mathcal{C}(ZM(m,n,r))$ and the lattice of divisors of $mn$, we have 
\begin{equation}\label{r1}
[H(m_i,n_i,0)]\cap [H(m_{j},n_{j},0)]=[H([m_i,m_j], [n_i, n_j], 0)].
\end{equation} 
We end this Section by noticing that the size of the conjugacy class $[H(m_i,n_i,0)]\in \mathcal{C}(ZM(m,n,r))$ is given by 
\begin{equation}\label{r2}
|[H(m_i,n_i,0)]|=\bigg(m_i, \frac{r^n-1}{r^{n_i}-1}\bigg ),
\end{equation}
as it was proved in \cite{1}.   

\section{The (cyclic) factorization number of a ZM-group}

Our main purpose in this Section is to find some formulas which allow us to compute explicitly the (cyclic) factorization number of a ZM-group. We will obtain two expressions which involve the sizes of some conjugacy classes of a ZM-group. Hence, the computation of the (cyclic) factorization number is reduced to summing some products of greatest common divisors of two positive integers.\\ 

\textbf{Theorem 3.1.} \textit{The factorization number of $ZM(m,n,r)$ is $$F_2(ZM(m,n,r))=\sum\limits_{\substack {m_2,m_3|m \\ n_2,n_3|n \\ (m_2,m_3)=1 \\ (n_2,n_3)=1}}\bigg( m_2,\frac{r^n-1}{r^{n_2}-1}\bigg)\bigg( m_3,\frac{r^n-1}{r^{n_3}-1}\bigg).$$
In particular, we have $$F_2(ZM(9,4,8))=85.$$}

\textbf{Proof.} Using the correspondence between the subgroup lattice $L(Z(m,n,r))$ and the set $L$, we can identify the entire group $Z(m,n,r)$ with $H(1,1,0)$. We are interested in counting the pairs $(H(m_2,n_2,x), H(m_3,n_3,y))\in L(H(1,1,0))^2$ such that $H(m_2,n_2,x)H(m_3,n_3,y)=H(1,1,0)$. We have $$H(1,1,0)=H(m_2,n_2,x)H(m_3,n_3,y)\Leftrightarrow |H(1,1,0)|=|H(m_2,n_2,x)H(m_3,n_3,y)|$$ $$\Leftrightarrow mn=\frac{\frac{mn}{m_2n_2}\frac{mn}{m_3n_3}}{\frac{mn}{[m_2,m_3][n_2,n_3]}}.$$ The last equality holds according to (\ref{r1}). Moreover, since $(m,n)=1$, the relation is further equivalent to  
$$\left\{\begin{array}{ll}
\left[  m_2,m_3\right] =m_2 m_3  \\
\left[  n_2,n_3\right] =n_2 n_3 
\end{array} \right.\Leftrightarrow \left\{\begin{array}{ll}
\left(  m_2,m_3\right) =1  \\
\left(  n_2,n_3\right) =1 
\end{array} \right. $$  
We recall that any two conjugate subgroups are isomorphic and each conjugacy class of a ZM-group contains its subgroups that have the same order. Hence, if $(H(m_2,n_2,x), H(m_3,n_3,y))$ is a factorization of $H(1,1,0)$, then $(H_2,H_3)$ is also such a factorization for any $H_2\in [H(m_2,n_2,0)]$ and $H_3\in [H(m_3,n_3,0)]$. These facts lead us to expressing the factorization number of $ZM(m,n,r)$ as follows:
$$F_2(ZM(m,n,r))=|\lbrace (H_2,H_3)\in [H(m_2,n_2,0)]\times [H(m_3,n_3,0)] \ | \ (m_2,m_3)=1, \ (n_2,n_3)=1\rbrace|,$$
where $m_2,m_3$ and $n_2,n_3$ are divisors of $m$ and $n$, respectively. Using (\ref{r2}), we obtain the desired explicit formula since 
\begin{align*}
F_2(ZM(m,n,r)) &=\sum\limits_{\substack {m_2,m_3|m \\ n_2,n_3|n \\ (m_2,m_3)=1 \\ (n_2,n_3)=1}}|[H(m_2,n_2,0)]||[H(m_3,n_3,0)]| \\
&= \sum\limits_{\substack {m_2,m_3|m \\ n_2,n_3|n \\ (m_2,m_3)=1 \\ (n_2,n_3)=1}}\bigg( m_2,\frac{r^n-1}{r^{n_2}-1}\bigg)\bigg( m_3,\frac{r^n-1}{r^{n_3}-1}\bigg).
\end{align*}
\hfill\rule{1,5mm}{1,5mm} \\
   
The same argument can be repeated to obtain the explicit formula for the cyclic factorization number of a ZM-group. This time we are working only with cyclic subgroups, so the divisors $m_2, m_3$ of $m$ and the divisors $n_3,n_3$ of $n$ must satisfy the supplementary property given by $L_1$, i.e. $\frac{m}{m_2}|r^{n_2}-1$ and $\frac{m}{m_3}|r^{n_3}-1$.\\

\textbf{Theorem 3.2.} \textit{The cyclic factorization number of $ZM(m,n,r)$ is $$CF_2(ZM(m,n,r))=\sum\limits_{\substack{m_2,m_3|m \\ n_2,n_3|n \\ (m_2,m_3)=1 \\ (n_2,n_3)=1 \\ \frac{m}{m_2}|r^{n_2}-1, \frac{m}{m_3}|r^{n_3}-1}}\bigg( m_2,\frac{r^n-1}{r^{n_2}-1}\bigg)\bigg( m_3,\frac{r^n-1}{r^{n_3}-1}\bigg).$$
In particular, we have $$CF_2(ZM(9,4,8))=36.$$}

\section{The (cyclic) subgroup commutativity degree of a ZM-group}

As we stated in Section 2, the subgroup commutativity degree of a finite group can be computed explicitly if we know all the factorization numbers of its subgroups. More exactly, we have 
$$sd(G)=\frac{1}{|L(G)|^2}\sum\limits_{H\in L(G)}F_2(H).$$
We mention that, for a ZM-group, the total number of its subgroups was found in \cite{1} and it is given by
$$|L(ZM(m,n,r))|=\sum\limits_{m_1|m}\sum\limits_{n_1|n}\bigg (m_1,\frac{r^{n}-1}{r^{n_1}-1}\bigg).$$ Hence, to obtain the subgroup commutativity degree of $Z(m,n,r)$  we must find an explicit form for the sum $\sum\limits_{H(m_1,n_1,s)\in L(ZM(m,n,r))}F_2(H(m_1,n_1,s)).$ Two conjugate subgroups have the same factorization number because they are isomorphic. Hence, we may write the same sum as $$\sum\limits_{[H(m_1,n_1,0)]\in \mathcal{C}(ZM(m,n,r))}|[H(m_1,n_1,0)]|F_2(H(m_1,n_1,0)).$$ We will denote by $f_{m_1,n_1}$ the quantity $F_2(H(m_1,n_1,0))$. Our next aim is to provide a formula which allows us to compute this parameter. We will denote by $(c_1)$, $(c_2)$, $(c_3)$ the conditions:

$$ (c_1) \left\{\begin{array}{ll}
m_2,m_3 |m  \\
n_2,n_3 |n \\
(m_2,m_3)=m_1 \\
(n_2,n_3)=n_1 \\
|[H(m_2,n_2,0)]|=1\\
|[H(m_3,n_3,0)]|=1,
\end{array} \right. (c_2) \left\{\begin{array}{ll}
m_2,m_3 |m  \\
n_2,n_3 |n \\
(m_2,m_3)=m_1 \\
(n_2,n_3)=n_1 \\
|[H(m_2,n_2,0)]|\neq 1\\
|[H(m_3,n_3,0)]|= 1,
\end{array} \right. (c_3) \left\{\begin{array}{ll}
m_2,m_3 |m  \\
n_2,n_3 |n \\
(m_2,m_3)=m_1 \\
(n_2,n_3)=n_1 \\
|[H(m_2,n_2,0)]|\neq 1\\
|[H(m_3,n_3,0)]|\neq 1.
\end{array} \right.$$ 
 
\textbf{Lemma 4.1.} \textit{Let $ZM(m,n,r)$ be a ZM-group and $H(m_1,n_1,0)$ be one of its subgroups. Then
\begin{align*}
f_{m_1,n_1}&=\sum\limits_{\substack (c_1)} \bigg( m_2, \frac{r^n-1}{r^{n_2}-1}\bigg)\bigg( m_3, \frac{r^n-1}{r^{n_3}-1}\bigg)+2\sum\limits_{\substack (c_2)}\frac{\big( m_2, \frac{r^n-1}{r^{n_2}-1}\big)\big( m_3, \frac{r^n-1}{r^{n_3}-1}\big)}{\big( m_1, \frac{r^n-1}{r^{n_1}-1}\big)}\\ &+\sum\limits_{\substack (c_3)}\frac{\big( m_2, \frac{r^n-1}{r^{n_2}-1}\big)\big( m_3, \frac{r^n-1}{r^{n_3}-1}\big)}{\big( m_1, \frac{r^n-1}{r^{n_1}-1}\big)^2}.
\end{align*}}  
\textbf{Proof.} 
We begin our proof by characterizing the fact that $(H(m_2,n_2,0), H(m_3,n_3,0))$ is a factorization of $H(m_1,n_1,0)$. Using the isomorphism between the lattice of divisors of $mn$ and $\mathcal{C}(ZM(m,n,r))$ and taking into account that $(m,n)=1$, we have
$$H(m_1,n_1,0)=H(m_2,n_2,0)H(m_3,n_3,0)\Leftrightarrow \left\{\begin{array}{ll}
H(m_2,n_2,0), H(m_3,n_3,0)\subseteq H(m_1,n_1,0)  \\
|H(m_1,n_1,0)|=|H(m_2,n_2,0)H(m_3,n_3,0)| 
\end{array} \right.$$ $$\Leftrightarrow \left\{\begin{array}{ll}
m_1|m_2, m_1|m_3  \\
n_1|n_2, n_1|n_3 \\
(m_2,m_3)=m_1 \\
(n_2,n_3)=n_1 
\end{array} \right. \Leftrightarrow \left\{\begin{array}{ll}
(m_2,m_3)=m_1\\
 (n_2,n_3)=n_1 
\end{array} \right.$$
Since conjugate subgroups are isomorphic, if the pair $(H(m_2,n_2,0), H(m_3,n_3,0))$ is a factorization of  $H(m_1,n_1,0)$, then $(H_2,H_3)$ is a factorization of $H_1$, where $H_i\in [H(m_i,n_i,0)]$, for $i\in \lbrace 1,2,3\rbrace$. If $[H(m_2,n_2,0)]$ and $[H(m_3,n_3,0)]$ are trivial, then $(H(m_2,n_2,0), H(m_3,n_3,0))$ is a factorization for any $H_1\in [H(m_1,n_1,0)]$. But, for $i\in\lbrace 2,3\rbrace$, if $[H(m_i,n_i,0)]$ is not a trivial conjugacy class, there are $\frac{|[H(m_i,n_i,0)]|}{|[H(m_1,n_1,0)]|}$ conjugates contained in $[H(m_i,n_i,0)]$ which contribute to the number of factorizations of a conjugate $H_1\in [H(m_1,n_1,0)]$. Again, this fact is a consequence of the existent isomorphisms between the subgroups which form the conjugacy class $[H(m_1,n_1,0)].$ Taking into account these aspects and the properties of the divisors $m_1,m_2,m_3$ and $n_1,n_2,n_3$ of $m$ and $n$, respectively, we obtain the desired explicit form of $f_{m_1,n_1}.$  
\hfill\rule{1,5mm}{1,5mm} \\

Using \textbf{Lemma 4.1.}, the quantities $f_{m_1,n_1}$ can be computed for each divisors $m_1$ of $m$ and $n_1$ of $n$. Consequently, the subgroup commutativity degree of a ZM-group is given by the following result.\\

\textbf{Theorem 4.2.} \textit{The subgroup commmutativity degree of $ZM(m,n,r)$ is  
$$sd(ZM(m,n,r))=\frac{\sum\limits_{m_1|m}\sum\limits_{n_1|n}\big(m_1,\frac{r^n-1}{r^{n_1}-1}\big)f_{m_1,n_1}}{\bigg[\sum\limits_{m_1|m}\sum\limits_{n_1|n}\big( m_1,\frac{r^{n}-1}{r^{n_1}-1}\big)\bigg]^2}.$$
In particular, we have $$sd(ZM(9,4,8))=\frac{13}{19}.$$}

The number of cyclic subgroups of a ZM-group was obtained in \cite{18} and its value is 
$$|L_1(ZM(m,n,r))|=\sum\limits_{m_1|m}\sum\limits_{\substack{n_1|n \\ \frac{m}{m_1}|r^{n_1}-1}}\bigg( m_1,\frac{r^{n}-1}{r^{n_1}-1}\bigg).$$
The connection between the cyclic subgroup commutativity degree of $ZM(m,n,r)$ and the cyclic factorization numbers of its subgroups is given by 
$$csd(ZM(m,n,r))=\frac{1}{|L_1(ZM(m,n,r))|^2}\sum\limits_{[H(m_1,n_1,0)]\in \mathcal{C}(ZM(m,n,r))}|[H(m_1,n_1,0)]|CF_2(H(m_1,n_1,0)).$$
If we denote the quantity $CF_2(H(m_1,n_1,0))$ by $cf_{m_1,n_1}$ and we repeat the arguments which were used in the proof of \textbf{Lemma 4.1.}, we obtain the following two analogous results which allow us to compute explicitly the cyclic subgroup commutativity degree of a ZM-group.\\ 

\textbf{Lemma 4.3.} \textit{Let $ZM(m,n,r)$ be a ZM-group and $H(m_1,n_1,0)$ be one of its subgroups. Then
\begin{align*}
cf_{m_1,n_1}&= \sum\limits_{\substack {(c_1) \\ \frac{m}{m_2}|r^{n_2}-1\\ \frac{m}{m_3}|r^{n_3}-1 \\ }} \bigg( m_2, \frac{r^n-1}{r^{n_2}-1}\bigg)\bigg( m_3, \frac{r^n-1}{r^{n_3}-1}\bigg)+2\sum\limits_{\substack {(c_2) \\ \frac{m}{m_2}|r^{n_2}-1\\ \frac{m}{m_3}|r^{n_3}-1 }}\frac{\big( m_2, \frac{r^n-1}{r^{n_2}-1}\big)\big( m_3, \frac{r^n-1}{r^{n_3}-1}\big)}{\big( m_1, \frac{r^n-1}{r^{n_1}-1}\big)}\\ &+\sum\limits_{\substack {(c_3)\\ \frac{m}{m_2}|r^{n_2}-1\\ \frac{m}{m_3}|r^{n_3}-1}}\frac{\big( m_2, \frac{r^n-1}{r^{n_2}-1}\big)\big( m_3, \frac{r^n-1}{r^{n_3}-1}\big)}{\big( m_1, \frac{r^n-1}{r^{n_1}-1}\big)^2}.
\end{align*}}

\textbf{Theorem 4.4.} \textit{The cyclic subgroup commutativity degree of $ZM(m,n,r)$ is 
$$csd(ZM(m,n,r))=\frac{\sum\limits_{m_1|m}\sum\limits_{n_1|n}\big(m_1,\frac{r^n-1}{r^{n_1}-1}\big)cf_{m_1,n_1}}{\bigg[\sum\limits_{m_1|m}\sum\limits_{\substack{n_1|n \\ \frac{m}{m_1}|r^{n_1}-1}}\big( m_1,\frac{r^{n}-1}{r^{n_1}-1}\big)\bigg]^2}.$$
In particular, we have $$csd(ZM(9,4,8))=\frac{17}{25}.$$}

Before we prove our next result, we introduce the function $g:\mathbb{N}\longrightarrow \mathbb{N}$ given by $$g(n)=n\sum\limits_{\substack{r|n\\ s|n}}\frac{1}{(r,s)},$$ which was used in \cite{16} to express the subgroup commutativity degree of the dihedral group $D_{2m}$, where $m\geq 3$ is a positive integer. This function is multiplicative and for any prime number $p$ and $\alpha\in \mathbb{N}^{*}$, we have
$$g(p^{\alpha})=\frac{(2\alpha +1)p^{\alpha+2}-(2\alpha+3)p^{\alpha+1}+p+1}{(p-1)^2}.$$ 
Now, we will find an explicit form of the results given by \textbf{Theorem 4.2.} and \textbf{Theorem 4.4.}  for a particular class of ZM-groups.\\

\textbf{Proposition 4.5.} \textit{Let $ZM(m,n,r)$ be a ZM-group such that $n$ is a prime number. Then\\
\textit{(i)}  $$sd(ZM(m,n,r))=\frac{\tau(m)^2+2\tau(m)\sigma(m)+g(m)}{(\tau(m)+\sigma(m))^2},$$\\
\textit{(ii)} $$csd(ZM(m,n,r))=\frac{\tau(m)(\tau(m)+m)+m(\tau(m)+1)}{(\tau(m)+m)^2}.$$}

\textbf{Proof.} We remark that the properties $r^n\equiv 1 \ (mod \ m)$ and $(m,r-1)=1$ leads us to stating that $(m_1,\frac{r^n-1}{r-1})=m_1$, for any divisor $m_1$ of $m$.\\
\textit{(i)} Since $n$ is a prime  number, the subgroup commutativity degree of $ZM(m,n,r)$ is 
\begin{equation}\label{r3}
sd(ZM(m,n,r))=\frac{\sum\limits_{m_1|m}m_1f_{m_1,1}+\sum\limits_{m_1|m}f_{m_1,n}}{\bigg(\sum\limits_{m_1|m}m_1+\sum\limits_{m_1|m}1\bigg)^2}=\frac{\sum\limits_{m_1|m}m_1f_{m_1,1}+\sum\limits_{m_1|m}f_{m_1,n}}{(\sigma(m)+\tau(m))^2}.
\end{equation}
According to \textbf{Lemma 4.1.}, the quantities $f_{m_1,1}$ and $f_{m_1,n}$ are
$$f_{m_1,1}=2\sum\limits_{\substack {m_2,m_3|m  \\ (m_2,m_3)=m_1}}\frac{m_2}{m_1}+\sum\limits_{\substack {m_2,m_3|m \\ (m_2,m_3)=m_1}}\frac{m_2m_3}{m_1^{2}},$$
and
$$f_{m_1,n}=\sum\limits_{\substack{m_2,m_3|m \\ (m_2,m_3)=m_1}}1,$$
respectively. The right-hand side of (\ref{r3}) can be further reduced since 
\begin{align*}
\sum\limits_{m_1|m}m_1f_{m_1,1}&=2\sum\limits_{m_1|m}\sum\limits_{\substack{m_2,m_3|m \\ (m_2,m_3)=m_1}}m_2+\sum\limits_{m_1|m}\sum\limits_{\substack{m_2,m_3|m \\ (m_2,m_3)=m_1}}\frac{m_2m_3}{m_1} \\
&= 2\sum\limits_{m_2,m_3|m}m_2+\sum\limits_{m_2,m_3|m}[m_2,m_3]=2\tau(m)\sigma(m)+\sum\limits_{u,v|m}\bigg[\frac{m}{u},\frac{m}{v}\bigg]\\
&=2\tau(m)\sigma(m)+\sum\limits_{u,v|m}\frac{m}{(u,v)}=2\tau(m)\sigma(m)+g(m),  
\end{align*}
and $$\sum\limits_{m_1|m}f_{m_1,n}=\sum\limits_{m_1|m}\sum\limits_{\substack {m_2,m_3|m \\ (m_2,m_3)=m_1}}1=\sum\limits_{m_2,m_3|m}1=\tau(m)^2.$$
The proof is finished once the replacements are made in (\ref{r3}).\\
\textit{(ii)} The cyclic subgroup commutativity degree of $ZM(m,n,r)$ is given by 
$$csd(ZM(m,n,r))=\frac{\sum\limits_{m_1|m}m_1cf_{m_1,1}+\sum\limits_{m_1|m}cf_{m_1,n}}{\bigg(\sum\limits_{m_1|m}\sum\limits_{\frac{m}{m_1}|r-1}m_1+\sum\limits_{m_1|m}\sum\limits_{\frac{m}{m_1}|r^n-1}1\bigg)^2}=\frac{\sum\limits_{m_1|m}m_1cf_{m_1,1}+\sum\limits_{m_1|m}cf_{m_1,n}}{(m+\tau(m))^2},$$ since the only divisor of $m$ which divides $r-1$ is 1 and all divisors of $m$ divide the quantity $r^n-1$. Using \textbf{Lemma 4.3.}, we have
$$cf_{m_1,1}=2\sum\limits_{\substack {m_2,m_3|m  \\ (m_2,m_3)=m_1 \\ \frac{m}{m_2}|r-1,\frac{m}{m_3}|r^{n}-1}}\frac{m_2}{m_1}+\sum\limits_{\substack {m_2,m_3|m \\ (m_2,m_3)=m_1 \\ \frac{m}{m_2}|r-1,\frac{m}{m_3}|r-1}}\frac{m_2m_3}{m_1^{2}},$$
and
$$cf_{m_1,n}=\sum\limits_{\substack{m_2,m_3|m \\ (m_2,m_3)=m_1\\ \frac{m}{m_2}|r^n-1,\frac{m}{m_3}|r^n-1}}1.$$
Consequently, 
\begin{align*}
\sum\limits_{m_1|m}m_1cf_{m_1,1}&=2\sum\limits_{m_1|m}\sum\limits_{\substack {m_2,m_3|m  \\ (m_2,m_3)=m_1 \\ \frac{m}{m_2}|r-1,\frac{m}{m_3}|r^{n}-1}}m_2+\sum\limits_{m_1|m}\sum\limits_{\substack {m_2,m_3|m \\ (m_2,m_3)=m_1 \\ \frac{m}{m_2}|r-1,\frac{m}{m_3}|r-1}}\frac{m_2m_3}{m_1}\\ &=2m\sum\limits_{m_1|m}1+m=2m\tau(m)+m,
\end{align*}
and
$$\sum\limits_{m_1|m}cf_{m_1,n}=\sum\limits_{m_1|m}\sum\limits_{\substack{m_2,m_3|m \\ (m_2,m_3)=m_1\\ \frac{m}{m_2}|r^n-1,\frac{m}{m_3}|r^n-1}}1=\sum\limits_{m_2,m_3|m}1=\tau(m)^2.$$
Then, the cyclic subgroup commutativity degree of $ZM(m,n,r)$ is 
$$csd(ZM(m,n,r))=\frac{2m\tau(m)+m+\tau(m)^2}{(m+\tau(m))^2}=\frac{\tau(m)(\tau(m)+m)+m(\tau(m)+1)}{(\tau(m)+m)^2},$$ as desired. 
\hfill\rule{1,5mm}{1,5mm}\\

It is interesting that for a prime number $n$, the explicit formulas for the (cyclic) subgroup commutativity degree of $ZM(m,n,r)$ depend only on the parameter $m$. Therefore, if we consider another prime number $n_1$ and a corresponding value $r_1$ satisfying $$(m,n_1)=(m,r_1-1)=1 \text{ and } r_1^{n_1} \equiv 1 \ (mod \ m),$$ then the equalities $sd(ZM(m,n,r))=sd(ZM(m,n_1,r_1)),csd(ZM(m,n,r))=csd(ZM(m,n_1,r_1))$ hold. Consequently, there is an infinite number of ZM-groups having the same (cyclic) subgroup commutativity degree. 

We remark that $sd(ZM(m,2,m-1))=sd(D_{2m})$ and $csd(ZM(m,2,m-1))=csd(D_{2m})$ for any odd number $m\geq 3$. We expected this since $ZM(m,2,m-1)=D_{2m}$, for any $m\geq 3$ such that $m\equiv 1 \ (mod \ 2)$. Hence, our work constitute a generalization of the results which were proved in \cite{16,19} for dihedral groups. Also, if we still consider $n$ to be a prime, by applying similar techniques to the ones that were used in the proof of \textbf{Proposition 4.5.}, one can show that the results concerning the (cyclic) factorization number of $ZM(m,n,r)$, which were obtained in Section 3, generalize the formulas provided by \cite{11,20} for the same concepts associated to dihedral groups.

Our final result points out a new class of groups whose (cyclic) subgroup commutativity degree tends to $0$ as the orders of the groups tend to infinity.\\

\textbf{Corollary 4.6.} \textit{Let $n$ and $p\geq 3$ be two prime numbers and let $\alpha$ be a positive integer. Then\\ 
\textit{(i)}  $$\displaystyle \lim_{\alpha \to\infty}sd(ZM(p^{\alpha},n,r))=0,$$\\
\textit{(ii)} $$\displaystyle \lim_{\alpha \to\infty}csd(ZM(p^{\alpha},n,r))=0.$$}

\textbf{Proof.} Both results easily follow using the explicit formulas given by \textbf{Proposition 4.5.} in the particular case $m=p^{\alpha}$. We have $\tau(p^{\alpha})=\alpha+1, \sigma(p^{\alpha})=\frac{p^{\alpha+1}-1}{p-1}$ and we recall that 
$$g(p^{\alpha})=\frac{(2\alpha +1)p^{\alpha+2}-(2\alpha+3)p^{\alpha+1}+p+1}{(p-1)^2}.$$
Therefore,
$$\displaystyle \lim_{\alpha \to\infty}sd(ZM(p^{\alpha},n,r))=\displaystyle \lim_{\alpha \to\infty}\frac{\tau(p^{\alpha})^2+2\tau(p^{\alpha})\sigma(p^{\alpha})+g(p^{\alpha})}{(\tau(p^{\alpha})+\sigma(p^{\alpha}))^2}$$ 
$$=\displaystyle \lim_{\alpha \to\infty} \frac{[(4\alpha+3)(p-1)-2]p^{\alpha+1}+(\alpha+1)(p-1)[(\alpha+1)(p-1)-2]+p+1}{[p^{\alpha+1}+(\alpha+1)(p-1)-1]^2}=0$$
and
\begin{align*}
\displaystyle \lim_{\alpha \to\infty}csd(ZM(p^{\alpha},n,r))&=\displaystyle \lim_{\alpha \to\infty} \frac{\tau(p^{\alpha})(\tau(p^{\alpha})+p^{\alpha})+p^{\alpha}(\tau(p^{\alpha})+1)}{(\tau(p^{\alpha})+p^{\alpha})^2}\\ &=\displaystyle \lim_{\alpha \to\infty}\frac{(2\alpha+3)p^{\alpha}+(\alpha+1)^2}{(p^{\alpha}+\alpha+1)^2}=0.
\end{align*}
\hfill\rule{1,5mm}{1,5mm}

      
\section{Conclusions and further research}

The probabilistic aspects of finite groups theory constitute a fruitful research topic and the connections between them may establish new relevant results. We end our paper by indicating the following three open problems.\\

\textbf{Problem 5.1.} Consider the sets $\lbrace sd(G) \ | \ G=finite \ group\rbrace$ and $\lbrace csd(G) \ | \ G=finite \ group\rbrace$. Study their density in $[0,1]$.\\

\textbf{Problem 5.2.} For a fixed $\alpha\in (0,1)$, describe the structure of finite groups $G$ satisfying $sd(G)=(\leqslant, \geqslant)\alpha$ or $csd(G)=(\leqslant, \geqslant)\alpha$.\\

\textbf{Problem 5.3.} Do there exist finite groups $G$ such that $sd(G)=csd(G)\ne 1?$

\vspace*{3ex}
\small
\hfill
\begin{minipage}[t]{6cm}
Mihai-Silviu Lazorec \\
Faculty of  Mathematics \\
"Al.I. Cuza" University \\
Ia\c si, Romania \\
e-mail: {\tt Mihai.Lazorec@student.uaic.ro}
\end{minipage}

\end{document}